\documentclass{article}

\usepackage{amssymb}
\usepackage[pdftex]{graphicx}
\usepackage{color}

\newtheorem{theorem}{Theorem}[section]

\newtheorem{lemma}{Lemma}[section]

\title{Uniform Model Completeness for the Real Field with the Weierstrass $\wp$ Function}
\author{Ricardo Bianconi\thanks{Departamento de Matem\'atica,
Instituto de Matem\'atica e Estat\'istica da Universidade de S\~ao Paulo,
Rua do Mat\~ao, 1010, Cidade Universit\'aria, CEP 05508-090, S\~ao Paulo, SP, Brazil.
e-mail: bianconi@ime.usp.br
}
}
\date{}

\begin{document}

\maketitle

\begin{abstract}
In this work is we prove model completeness for the expansion of the real field by the Weierstrass $\wp$ function as a function of the variable $z$ and the parameter (or period) $\tau$. We need to existentially define the partial derivatives of the $\wp$ function with respect to the variable $z$ and the parameter $\tau$. In order to obtain this result we need to include in the structure function symbols for the unrestricted exponential function and restricted sine function, the Weierstrass $\zeta$ function and the quasimodular form $E_2$. We prove some auxiliary model completeness results with the same functions composed with appropriate change of variables. In the conclusion we make some remarks about the noneffectiveness of our proof and the difficulties to be overcome to obtain an effective model completeness result.

\textbf{Mathematical Subject Classification.} Primary: 03C10. Secondary: 03C64, 03C98, 14H52, 33E05.

\textbf{Keywords:} model completeness, Weierstrass systems, elliptic functions, $\wp$ function, Weierstrass $\zeta$ function, modular forms, o-minimality, definability.
\end{abstract}

\textbf{Short title:} {Uniform Model Completeness for the $\wp$ Function}

\section{Introduction}

O-minimal techniques have been used to prove important results in algebraic geometry. Johnathan Pila \cite{pila2011}  has proved in 2011 the Andr\'e-Oort Conjecture for $\mathbb{C}\sp{n}$, where he makes essential use of his work with Alex Wilkie \cite{pila-wilkie2006} on rational points on definable sets in an o-minimal structure, and work by Ya'acov Peterzil and Sergei Starchenko \cite{peterzil-starchenko-wp2004} on the uniform definability of Weierstrass' $\wp$ function in the o-minimal structure $\mathbb{R}_{\mathit{an}}$ (the expansion of the real field by real analytic functions restricted to $[-1,1]\sp{n}$). This structure is model complete, so the $\wp$ function is existentially definable. This suggests the question of which ``minimal'' reduct of $\mathbb{R}_{\mathit{an}}$ is still model complete and defines $\wp$.

We prove here model completeness of expansions of the real fields by $\wp$ function as a function of the variable $z$ and of the parameter $\tau$. The structure obtained is interpretable in $\mathbb{R}_{\mathit{an}}$ and so it is o-minimal. We present a model complete reduct of $\mathbb{R}_{\mathit{an,exp}}$ (with finitely many function symbols, and with the unrestricted exponential function) in which the function symbols of that structure are existentially interpretable (and conversely). This gives the desired result.
Because of the techniques we use, we need the exponential function, restricted sine function and modular forms.

We conjecture that the exponential function and the restricted sine function are necessary for such definition, that is, we conjecture that these functions are not \emph{existentially} definable from $\wp$ and $\wp'$ alone (although they are definable from $\wp$, see \cite[Theorem 5.7, p. 545]{peterzil-starchenko-wp2004} for the exponential function and \cite[Section 20.222, pp. 438-439]{whittaker-watson} for the sine function).

This work is organised as follows. In the following section we present the general techniques to prove model completeness for reducts of $\mathbb{R}_{\mathit{an,exp}}$. The main results are stated in the next section. The proofs are given in full in the sequel. In the concluding remarks we deal with the problem of finding an effective proof of these model completeness results.

\medskip

\textbf{Notation:} $\mathbb{Z}$ denotes the ring of rational integers, $\mathbb{N}$ the set of nonnegative integers, $\mathbb{R}$ the field of real numbers, $\mathbb{C}$ the field of complex numbers; $\mathfrak{Re}(z)$ and $\mathfrak{Im}(z)$ denote the real and imaginary parts of $z$, the letter $i$ denotes $\sqrt{-1}$,  and $\bar{x}$ denotes the tuple $(x_1,\dots,x_n)$ for some unspecified $n$. Other notations are explained in the text.

\section{Model Completeness Criteria} %%%% MODEL COMPLETENESS CRITERIA

We use the following criteria for model completeness. They are based on the reduction to the so-called Weierstrass Systems of the proof of model completeness of the structure $\mathbb{R}_{\mathit{an}}$ by Denef and van de Dries in \cite{denef-vandendries1988}, as suggested by van den Dries in \cite{vandendries1988} and presented in \cite{bianconi1991}. The versions below are based on these papers and in \cite{vandendries-etc-1996}, as presented in \cite{bianconi2014}.

We stress here the point that these results rely on the (topological) compactness of closed and bounded polyintervals in $\mathbb{R}\sp{n}$ and polydisks in $\mathbb{C}\sp{n}$, which gives a noneffective feature to the proofs. We discuss
this in the concluding remarks.

\medskip

Recall that the theory $T$ of a structure $M$ is called \textit{strongly model complete} if for each formula $\varphi(\bar{x})$ there is a quantifier free formula $\psi(\bar{x},\bar{y})$ such that
$$
T\vdash\forall \bar{x}\,[\varphi(\bar{x})\leftrightarrow\exists\bar{y}\,\psi(\bar{x},\bar{y})]
$$
and $T\vdash\forall\bar{x}\forall\bar{y}\forall\bar{z}(\psi(\bar{x},\bar{y})\wedge\psi(\bar{x},\bar{z})\to\bar{y}=\bar{z})$. We say that such formula is a \textit{strong (existential) formula} and the set it defines, \textit{strongly definable}.

The first criterion deals with functions restricted to compact polyintervals and is the main step to obtain the more general result.

\begin{theorem}[{\cite[Theorem 2]{bianconi2014}}]\label{model-completeness-W-system} 
Let $\hat{R}=\langle\mathbb{R},\mathit{constants},+,-,\cdot,<,(F_{\lambda})_{\lambda\in\Lambda}\rangle$ be an expansion of the field of real numbers, where for each $\lambda\in\Lambda$, $F_{\lambda}$ is the restriction to a compact polyinterval $D_{\lambda}\subseteq\mathbb{R}\sp{n_{\lambda}}$ of a real analytic function whose domain contains $D_{\lambda}$, and defined as zero outside $D_{\lambda}$, such that there exists a complex analytic function $g_{\lambda}$ defined in a neighbourhood of a polydisk $\Delta_{\lambda}\supseteq D_{\lambda}$ and such that
\begin{enumerate}
\item
$g_{\lambda}$ is strongly definable in $\hat{R}$ and the restriction of $g_{\lambda}$ to $D_{\lambda}$ coincides with $F_{\lambda}$ restricted to the same set;
\item
for each $a\in\Delta_{\lambda}$ there exists a compact polydisk $\Delta$ centred at $a$ and contained in the domain of $g_{\lambda}$, such that all the partial derivatives of the restriction of $g_{\lambda}$ to $\Delta$ are strongly definable in $\hat{R}$.
\end{enumerate}
Under these hypotheses, the theory of $\hat{R}$ is strongly model complete.
\end{theorem}

Now we introduce the unrestricted exponential function.

\begin{theorem}[{\cite[Theorem 4]{bianconi2014}}]\label{model-completeness-W-system-exp} 
Let $\hat{R}$ be the structure described in Theorem \ref{model-completeness-W-system}. We assume that the functions
$$
\exp\lceil_{[0,1]}(x)=
\left\{
\begin{array}{lcl}
\exp x & \mbox{if} & 0\leq x\leq 1,\\
0 & & \mbox{otherwise};
\end{array}
\right.
$$
$$
\sin\lceil_{[0,\pi]}(x)=
\left\{
\begin{array}{lcl}
\sin x & \mbox{if} & 0\leq x\leq \pi,\\
0 & & \mbox{otherwise},
\end{array}
\right.
$$
have representing function symbols in its language.
The expansion $\hat{R}_{\mathrm{exp}}$ of $\hat{R}$ by the inclusion of the (unrestricted) exponential function ``$\exp$'' is strongly model complete.
\end{theorem}

We use this theorem to prove model completeness of the expansion of the field of real numbers by functions defined in unbounded sets of the real numbers using the exponential function and the restricted sine function to existentially interpret that structure into another which satisfies its hypotheses.

%%%%%%% THE MAIN RESULTS %%%%%%%

\section{The Main Results}

In order to state the main theorems we need to introduce some notations and definitions. We follow closely the work \cite{peterzil-starchenko-wp2004} by Peterzil and Starchenko on the uniform definability of the $\wp$ functions on $\mathbb{R}_{\mathit{an}}$ and use their notation for easy comparison with their work.

We define $\mathfrak{F}_d=\{\tau\in\mathbb{H}:-1/2\leq\mathfrak{Re}(\tau)<1/2$, and $\mathfrak{Im}(\tau)\geq d\}$, where $d=\sqrt{3}/2$. Observe that this set contains the usual fundamental domain of the action of $SL(2,\mathbb{Z})$ on $\mathbb{H}$, the set $\mathfrak{F}=\{\tau\in\mathbb{H}:-1/2\leq\mathfrak{Re}(\tau)<1/2$, $|\tau|\geq 1\}$. For each $\tau\in\mathbb{H}$ we define the set $E_{\tau}=\{z\in\mathbb{C}:z=a\tau+b$, for $a,b\in\mathbb{R}$, $a,b\geq 0$ and $a,b< 1\}$ and the set $E\sp{\mathfrak{F}_d}=\{(\tau,z)\in\mathbb{H}\times\mathbb{C}:\tau\in\mathfrak{F}_d$ and $z\in E_{\tau}\}$. This set is semialgebraic and contains a fundamental region for the actions of $SL(2,\mathbb{Z})$ on $\mathbb{H}$ and the action of $\mathbb{Z}\times\mathbb{Z}$ on $\mathbb{C}$, generated by the translations $z\mapsto z+1$ and $z\mapsto z+\tau$.

We write $z=x+iy$ and $\tau=u+iv$, with $x,y,u,v$ real variables, and define the functions

$$
\mathit{RP}(u,v;x,y)=\left\{\begin{array}{lll}
\mathfrak{Re}(\wp)(\tau;z) & \mbox{if} &  (\tau,z)\in E\sp{\mathfrak{F}_d}\\
0 & & \mbox{otherwise},
\end{array}
\right.
$$
$$
\mathit{IP}(u,v;x,y)=\left\{\begin{array}{lll}
\mathfrak{Im}(\wp)(\tau;z) & \mbox{if} &  (\tau,z)\in E\sp{\mathfrak{F}_d}\\
0 & & \mbox{otherwise},
\end{array}
\right.
$$
the real and imaginary parts of Weierstrass' $\wp$ function,
$$
\mathit{RZ}(u,v;x,y)=\left\{\begin{array}{lll}
\mathfrak{Re}(\zeta)(\tau;z) & \mbox{if} &  (\tau,z)\in E\sp{\mathfrak{F}_d}\\
0 & & \mbox{otherwise},
\end{array}
\right.
$$
$$
\mathit{IZ}(u,v;x,y)=\left\{\begin{array}{lll}
\mathfrak{Im}(\zeta)(\tau;z) & \mbox{if} &  (\tau,z)\in E\sp{\mathfrak{F}_d}\\
0 & & \mbox{otherwise},
\end{array}
\right.
$$
the real and imaginary parts of Weierstrass' zeta function, which can be defined by
$$
\frac{\partial\zeta}{\partial z}(\tau;z)=-\wp(\tau;z),\ \mbox{and}\ \lim_{z\to 0}\left(\zeta(\tau;z)-\frac{1}{z}\right)=0;
$$
and also
$$
\mathit{RE}_2(u,v)=\left\{\begin{array}{lll}
\mathfrak{Re}(\tilde{E}_2)(\tau) & \mbox{if} &  \tau\in \mathfrak{F}_d\\
0 & & \mbox{otherwise},
\end{array}
\right.
$$
$$
\mathit{IE}_2(u,v)=\left\{\begin{array}{lll}
\mathfrak{Im}(\tilde{E}_2)(\tau) & \mbox{if} &  \tau\in \mathfrak{F}_d\\
0 & & \mbox{otherwise},
\end{array}
\right.
$$
the real and imaginary parts of the quasimodular form $E_2(\tau)$ (see the proof of Lemma \ref{lemma1} below for its precise definition).

We can state now the main result, where $\exp$ is the unrestricted exponential function and $\sin\lceil_{[0,\pi]}$ is the sine function restricted to the interval $[0,\pi]$ and defined as 0 outside this interval.

\begin{theorem}\label{main-thm}
The structure $\mathcal{R}_{\wp}=\langle \mathbb{R}$, \textrm{constants}, $+$, $-$, $\cdot$, $<$, $\mathit{RP}$, $\mathit{IP}$, $\mathit{RZ}$, $\mathit{IZ}$, $\mathit{RE}_2$, $\mathit{IE}_2$, $\exp$, $\sin\lceil_{[0,\pi]}\rangle$ is model complete.
\end{theorem}

The idea of the proof is to use the change of variables $q=\exp(2\pi i \tau)$, $u=\exp(2\pi iz)$, in order to interpret the structure $\mathcal{R}_{\wp}$ into an auxiliary structure which satisfies the hypotheses of Theorem \ref{model-completeness-W-system-exp}. This is done in the following section.

%%%%%%%%%

\section{Proof of Theorem \ref{main-thm}}

We divide this section into three steps, focussing in each item of the model completeness criteria.

\medskip

\textsc{Step 1.} The model completeness criteria require that all the partial derivatives of the strongly definable complex functions are also strongly definable. The following two lemmas show this.

\begin{lemma}\label{lemma1}
The field $\mathbb{C}(z,\tau,\wp,\wp',g_2,g_3,\zeta,E_2)$ is closed under differentiation.
\end{lemma}

\textbf{Proof:} We recall the differential equation satisfied by $\wp$, where $\wp'$ indicates the (partial) derivative with respect to $z$.
$$
\wp\sp{\prime\,2}(\tau;z)=4\wp(\tau;z)\sp{3}-g_2(\tau)\wp(\tau;z)-g_3(\tau),
$$
from which we can obtain all partiall derivatives with respect to $z$ as rational functions of $\wp$, $\wp'$, $g_2$ and $g_3$.

By definition, $\partial\zeta(\tau;z)/\partial z=-\wp(\tau;z)$, from which we can obtain all partiall derivatives of $\zeta$ with respect to the variable $z$.

There remains to differentiate with respect to the variable (or parameter) $\tau$.

We start with the functions $g_2$ and $g_3$, which are modular forms of degree 4 and 6, respectively. Here we follow \cite[Sections 2.1, 2.3, 5.1]{zagier2008}.

Write $g_2=60E_4$ and $g_3=140E_6$, where
$$
E_{2k}(\tau)=\sum_{(m,n)\in\mathbb{Z}\sp{2}\setminus\{(0,0)\}} \frac{1}{(m+n\tau)\sp{2k}},
$$
for $k=2,3$.

For $k=1$ the series does not converge absolutely but we can define
$$
G_2(\tau)=\zeta(2)E_2(\tau)=\frac{1}{2}\sum_{n\neq 0}\frac{1}{n\sp{2}}+\frac{1}{2}\sum_{m\neq 0}\sum_{n\in\mathbb{Z}}\frac{1}{(m\tau+n)\sp{2}},
$$
where the sums cannot be interchanged (the series converges conditionally). This defines a \textbf{quasimodular form of weight 2}, which has the transformation formula
$$
G_2\left(\frac{a\tau+b}{c\tau+d}\right)=(c\tau+d)\sp{2}G_2(\tau)-\pi ic(c\tau+d),
$$
which is not a modular form. Its importance is that the ring $\mathbb{C}[E_2,E_4,E_6]$ contains all the modular forms and is closed under differentiation because of Ramanujan's formulas
$$
E_2'=\frac{E_2\sp{2}-E_4}{12},\ E_4'=\frac{E_2E_4-E_6}{3},\ E_6'=\frac{E_2E_6-E_4\sp{2}}{2},
$$
where $E'_k(\tau)$ denotes $(1/2\pi i)dE_k/d\tau$ (see \cite[Proposition 15, p. 49]{zagier2008}).

For the functions $\wp$, $\wp'$ and $\zeta$ we use the following formulas \cite[18.6.19-22, p. 641]{abramowitz-stegun1972}, where we omit showing the variables of the functions $\wp$, $\wp'$ and $\zeta$, and consider $g_2$ and $g_3$ as variables and denote $\Delta=g_2\sp{3}-27g_3\sp{2}$.
$$
\Delta\frac{\partial\wp}{\partial g_3}=\left(3g_2\zeta-\frac{9}{2}g_3z\right)\wp'+6g_2\wp\sp{2}-9g_3\wp -g_2\sp{2}\phantom{mmx}
$$
$$
\Delta\frac{\partial\wp}{\partial g_2}=\left(-\frac{9}{2}g_3\zeta+\frac{g_2\sp{2}z}{4}\right)\wp'-9g_3\wp\sp{2}+\frac{g_2\sp{2}}{2}\wp+\frac{3}{2}g_2g_3\phantom{x}
$$
$$
\Delta\frac{\partial\zeta}{\partial g_3}=-3\zeta\left(g_2\wp+\frac{3g_3}{2}\right)+\frac{z}{2}\left(9g_3\wp+\frac{g_2\sp{2}}{2}\right)-\frac{3}{2}g_2\wp'
$$
$$
\Delta\frac{\partial\zeta}{\partial g_2}=\frac{\zeta}{2}\left(9g_3\wp+\frac{g_2\sp{2}}{2}\right)-\frac{g_2z}{2}\left(\frac{g_2\wp}{2}+\frac{3g_3}{4}\right)+\frac{9}{4}g_3\wp'
$$

From these we can obtain the derivatives of $\wp'$ with respect to $g_2$ and $g_3$. Because of these formulas we see the need of the $\zeta$ function.

The chain rule gives the derivatives of $\wp$, $\wp'$ and $\zeta$ with respect to the variable $\tau$, as required.\hspace{\fill}$\square$

\medskip

As an immediate consequence, we have

\begin{lemma}
The partial derivatives of all orders of the real and imaginary parts of the functions $\wp$ and $\zeta$, and form $E_2$ are strongly definable from their real and imaginary parts.
\end{lemma}

\textbf{Proof:} The only missing formulas are $g_2=-4(e_1e_2+e_1e_3+e_2e_3)$ and $g_3=4e_1e_2e_3$, where $e_1=\wp(\tau;1/2)$, $e_2=\wp(\tau;\tau/2)$ and $e_3=\wp(\tau;(1+\tau)/2)$ are the roots of the polynomial $4\wp\sp{3}-g_2\wp-g_3$. The correct choice of the branch of the square root of this polynomial provides the definability of $\wp'$.

If we write $z=x+iy$ and $\tau=\alpha+i\beta$, we have  $\wp(\alpha+i\beta;x+iy)=\mathit{RP}(\alpha,\beta;x,y)+i\mathit{IP}(\alpha,\beta;x,y)$, $\zeta=\mathit{RZ}(\alpha,\beta;x,y)+i\mathit{IZ}(\alpha,\beta;x,y)$ and $E_2(\tau)=\mathit{RE}_2(\alpha,\beta)+i \mathit{IE}_2(\alpha,\beta)$. From these formulas we can obtain the strong definability of all partial derivatives of $\wp$, $\zeta$ and $E_2$ with respect to both variables $z$ and $\tau$.\hspace{\fill}$\square$

\medskip

\textsc{Step 2.} The model completeness criteria also require that the functions are definably extended to a neighbourhood of the original domains.

\begin{lemma}
Extensions of the function $E_2$ can be strongly defined from $E_2$. Extensions of the functions $\wp$ and $\zeta$ in both variables can be strongly defined from these functions.
\end{lemma}

\textbf{Proof:}
Recall that $E_2(\tau\pm 1)=E_2(\tau)$ and, in general,
$$
G_2\left(\frac{a\tau+b}{c\tau+d}\right)=(c\tau+d)\sp{2}G_2(\tau)-\pi ic(c\tau+d),
$$
which allows the definition of extensions of $E_2$ beyond the fundamental domain.

The function $\wp$ satisfies the addition formulas
$$
\wp(\tau;u+v)=\frac{1}{4}\left(\frac{\wp'(\tau;u)-\wp'(\tau;v)}{\wp(\tau;u)-\wp(\tau;v)}\right)\sp{2}-\wp(\tau;u)-\wp(\tau;v),
$$
and, for the case where $u=v$,
$$
\wp(\tau;2u)=-2\wp(\tau;u)-\left(\frac{\wp''(\tau;u)}{2\wp'(\tau;u)}\right)\sp{2}.
$$
It  also satisfies $\wp(\tau+1;z)=\wp(\tau;z)$ and
$\wp(-1/\tau;z)=\tau\wp(\tau;\tau z)$.
We recall that the M\"obius transformations $S(\tau)=-1/\tau$ and $T(\tau)=\tau+1$ generate the group $SL(2,\mathbb{Z})$. With these equations we can strongly define the desired extensions of $\wp$.

The function $\zeta$ satisfies $\zeta(\tau;z+\omega)=\zeta(\tau;z)+\eta(\omega)$, where $\eta:\Lambda_{\tau}\to\mathbb{C}$ is  group a homomorphism from the period lattice $\Lambda_{\tau}=\{m+n\tau:m,n\in\mathbb{Z}\}$ into $\mathbb{C}$. For $\omega\in\Lambda_{\tau}\setminus2\Lambda_{\tau}$, we have the formula $\eta(\omega)=2\zeta(\tau;\omega/2)$ (see \cite[Proposition 5.2, pp. 40-41]{silverman1994}). The zeta function also satisfies a kind of summation formula \cite[Section 20.41, \textit{Example 1}, p. 446]{whittaker-watson}
$$
[\zeta(\tau;x)+\zeta(\tau;y)+\zeta(\tau;z)]\sp{2}+
\wp(\tau;x)+\wp(\tau;y)+\wp(\tau;z)=0,
$$
if $x+y+z=0$. From the series defining the zeta function we obtain $\zeta(\tau+1;z)=\zeta(\tau;z)$ and $\zeta(-1/\tau;z)=\tau\zeta(\tau;\tau z)$.

This proves the lemma.\hspace{\fill}$\square$

\medskip

%%%%%%%%% CHANGE OF VARIABLES

\textsc{Step 3.} The model completeness criteria require that the functions are defined in compact polydisks, with the exception of the exponential function.
We apply now the change of variables $q=\exp(2\pi i \tau)$, $u=\exp(2\pi iz)$, which transforms the unbounded domains into bounded ones. The lemmas above together with this transfomation of variables show that we can cover the transformed set by finitely many polydisks, where the relevant functions admit strongly existentially definable extensions and, therefore, we can apply Theorem \ref{model-completeness-W-system-exp}.

\medskip

Recall the sets $\mathfrak{F}_d=\{\tau\in\mathbb{H}:-1/2\leq\mathfrak{Re}(\tau)<1/2$, and $\mathfrak{Im}(\tau)\geq d\}$, $d=\sqrt{3}/2$;
for each $\tau\in\mathfrak{F}_d$,
$E_{\tau}=\{z\in\mathbb{C}:z=a\tau+b$, for $a,b\in\mathbb{R}$, $a,b\geq 0$ and $a,b< 1\}$; and
$E\sp{\mathfrak{F}_d}=\{(\tau,z)\in\mathbb{H}\times\mathbb{C}:\tau\in\mathfrak{F}_d$, and $z\in E_{\tau}\}$.

We define the sets $\tilde{E}_{\tau}=\{z\in\mathbb{C}:$ there are $0\leq a,b\leq 1$, with $a+b\leq 1$, such that  $z=1/8+a/4+\tau/8+b\tau/4\}$, the covex closure of the paralellogram with vertices $(1+\tau)/8$, $(3+\tau)/8$, $(1+3\tau)/8$ and $(3+3\tau)/8$. Define
 $\tilde{E}\sp{\mathfrak{F}_d}=\{(\tau,z)\in\mathbb{H}\times\mathbb{C}:\tau\in\mathfrak{F}_d$, and $z\in \tilde{E}_{\tau}\}$.

The image of the set $\tilde{E}\sp{\mathfrak{F}_d}$ under the mapping $(\tau;z)\mapsto (q;u)$ is the set $M_{\delta}=\{(q,u)\in\mathbb{C}:$ $0<|q|\leq\delta$ and $(\log u)/(2\pi i) \in \tilde{E}_{(\log q)/2\pi i} \}$, where $\delta=\exp(-2\pi d)<~1$.

We apply this change of variables to the functions $\wp$, $\zeta$ and $E_2$.

For the $\wp$-function we have \cite[Theorem 6.2, p.50]{silverman1994}
$$
\wp_0(q;u)=(2\pi i)\sp{2}\left[\sum_{m\in\mathbb{Z}}\frac{uq\sp{m}}{(1-uq\sp{m})\sp{2}}+\frac{1}{12}-2\sum_{m\geq 1}\frac{q\sp{m}}{(1-q\sp{m})\sp{2}}\right],
$$

For the zeta function we have \cite[Theorem 6.3, pp. 52-53]{silverman1994}
$$
\frac{1}{2\pi i}\zeta_0(q;u)=\sum_{n\geq 0}\frac{-q\sp{n}u}{1-q\sp{n}u}+\sum_{n\geq 1}\frac{q\sp{n}u\sp{-1}}{1-q\sp{n}u\sp{-1}}+\frac{1}{2\pi i}\eta(1)z-\frac{1}{2},
$$
where we have a logarithmic term $z=(\log u)/2\pi i$ and
$$
\frac{1}{(2\pi i)\sp{2}}\eta(1)=\frac{1}{12}\left[ -1+24\sum_{n\geq 1}\frac{q\sp{n}}{(1-q\sp{n})\sp{2}}\right]
$$

We see that these series converge in an open neighbourhood of the image $M_{\delta}$ and also converge if $q=0$. Observe that because for all $\tau\in\mathfrak{F}_d$ and $z\in \tilde{E}_{\tau}$, $z$ is bounded away from zero and so the logarithmic term in the series of $\zeta_0$ has no singularity in some open neighbourhood of $M_{\delta}$.

Write $q=q_0+iq_1$ and $u=u_0+iu_1$. We put $\delta=\exp(-2\pi d)$, with $d=\sqrt{3}/2$, and define the functions
$$
\mathit{RP}_0(q_0,q_1;u_0,u_1)=\left\{\begin{array}{lll}
\mathfrak{Re}(\wp_0)(q;u) & \mbox{if} &  (q_0,q_1;u_0,u_1)\in M_{\delta}\\
0 & & \mbox{otherwise},
\end{array}
\right.
$$
$$
\mathit{IP}_0(q_0,q_1;u_0,u_1)=\left\{\begin{array}{lll}
\mathfrak{Im}(\wp_0)(q;u) & \mbox{if} &  (q_0,q_1;u_0,u_1)\in M_{\delta}\\
0 & & \mbox{otherwise},
\end{array}
\right.
$$
the real and imaginary parts of Weierstrass' $\wp$ function,
$$
\mathit{RZ}_0(q_0,q_1;u_0,u_1)=\left\{\begin{array}{lll}
\mathfrak{Re}(\zeta_0)(q;u) & \mbox{if} &  (q_0,q_1;u_0,u_1)\in M_{\delta}\\
0 & & \mbox{otherwise},
\end{array}
\right.
$$
$$
\mathit{IZ}_0(q_0,q_1;u_0,u_1)=\left\{\begin{array}{lll}
\mathfrak{Im}(\zeta_0)(q;u) & \mbox{if} &  (q_0,q_1;u_0,u_1)\in M_{\delta}\\
0 & & \mbox{otherwise},
\end{array}
\right.
$$

We define the function $\tilde{E}_2(q)$ as the Fourier transform of the function $E_2(\tau)$, which is a power series in the variable $q$ converging in an open neighbourhood of the closed disk $\Delta_{\delta}$ centered at $q=0$ (see \cite[p. 19]{zagier2008}), and define
$$
\mathit{RE}_2\sp{*}(q_0,q_1)=\left\{\begin{array}{lll}
\mathfrak{Re}(\tilde{E}_2)(q) & \mbox{if} &  (q_0,q_1)\in \mathfrak{F}_d\\
0 & & \mbox{otherwise},
\end{array}
\right.
$$
$$
\mathit{IE}_2\sp{*}(q_0,q_1)=\left\{\begin{array}{lll}
\mathfrak{Im}(\tilde{E}_2)(q) & \mbox{if} &  (q_0,q_1)\in \mathfrak{F}_d\\
0 & & \mbox{otherwise},
\end{array}
\right.
$$

\medskip

We apply the model completeness criteria from Theorem \ref{model-completeness-W-system} and the previous lemmas to conclude

\begin{theorem}\label{main-aux}
The structure $\mathcal{R}_{\mathit{aux}}=\langle \mathbb{R}$, \textrm{constants}, $+$, $-$, $\cdot$, $<$, $\mathit{RP}_0$, $\mathit{IP}_0$, $\mathit{RZ}_0$, $\mathit{IZ}_0$, $\mathit{RE}_2\sp{*}$, $\mathit{IE}_2\sp{*}$, $\exp\lceil_{[0,1]}$, $\sin\lceil_{[0,\pi]}\rangle$ is model complete.\hspace{\fill}$\square$
\end{theorem}

\textbf{Proof of Theorem \ref{main-thm}:}

By the model completeness Theorem \ref{model-completeness-W-system-exp}, the expansion of the structure $\mathcal{R}_{\mathit{aux}}$ including the full exponential function is strongly model complete. With this result at hand it is easy to see that the structure $\mathcal{R}_{\wp}$ is also strongly model complete.

Indeed, a definable set in the structure $\mathcal{R}_{\wp}$ can be mapped onto a set definable in $\mathcal{R}_{\mathit{aux},\exp}$ and so, it can be (strongly) existentially definable. This can be mapped back onto the original set and preserving the strong definability, because the map and its inverse are both strongly definable.\hspace*{\fill}$\square$

\medskip

\textit{Remark.} If we take the limit of $\wp(\tau;z)$ for $\tau\to i\infty$, with $|\mathfrak{Re}(\tau)|\leq 1/2$, we obtain the series $$
\sum_{n=-\infty}\sp{+\infty}\frac{1}{(z-n)\sp{2}}-2\sum_{n=1}\sp{\infty}\frac{1}{n\sp{2}}=
\frac{\pi\sp{2}}{\sin\sp{2}(\pi z)}-\frac{\ \pi\sp{2}}{3}
$$
(see \cite[Section 20.222, pp. 438-439]{whittaker-watson} and \cite[Formula 4.3.92, p. 75]{abramowitz-stegun1972}). Therefore, unless the restricted sine function  were (strongly, or just existentially) definable from the functions $\wp$, $\zeta$, $E_2$ and $\exp$, the reduct of the structure $\mathcal{R}_{\wp}$ without the restricted sine function would not be model complete. Actually, we conjecture that without either the sine or the exponential functions the (reduct of the) structure $\mathcal{R}_{\wp}$ may not be model complete. Observe that both functions are definable from the function $\wp$ (\cite[Theorem 5.7, p. 545] {peterzil-starchenko-wp2004} shows the definability of $\exp$).
 
\section{Concluding Remarks}

Our method is based on the fact that the interval $[-1,1]$ is a compact set and so the proofs bear a noneffective feature.

Recently, Angus Macintyre  \cite{macintyre2003} has shown how to prove an effective version of the model completeness of the expansion of the field of the real numbers by the $\wp$ function with parameter $\tau=i$. His proof can be modified to prove effective model completeness at least for the cases where $\wp$ admits complex multiplication and, perhaps, for all cases.
He uses the inverse function of $\wp$, which is an elliptic integral and thus it is a Pfaffian function (see details in \cite{macintyre2008}), so it is amenable to the treatment used in his work with Alex Wilkie on the decidability of the real exponential field, \cite{macintyre-wilkie1996}. In the latter paper they prove an effective model completeness result for the expansion of the real field with the restricted exponential function, from which we can pinpoint three main ideas:
\begin{enumerate}
\item
Pfaffian functions are amenable to good induction proofs on their complexity;
\item
there exist computable bounds on the number of connected components of zero sets of Pfaffian functions;
\item
the structure is polynomially bounded and they showed that there is a computable estimate on the behaviour at infinity of the one-variable definable functions.
\end{enumerate}

Their methods do not seem to apply to our case, because the derivative of $\wp$ with respect to the parameter $\tau$ involves the modular forms $E_4$ and $E_6$ which are not Pfaffian functions (but are \textit{noetherian} functions) and satisfy nonlinear differential equations involving $E_2$, $E_4$ and $E_6$ in a (to our knowledge) nontriangularisable way.

Of the three items stated above, only the third one (behaviour at infinity) seems to be more easily achievable.
 There is some work towards the first and second item (computable bounds on the number of connected components) surveyed in \cite{gabrielov-vorobjov2004}.

\bigskip

\textsc{Acknowledgements}

The author obtained good insights for this work from conversations with his coleagues Oswaldo Rio Branco de Oliveira, who works on finding elementary proofs of theorems in complex analysis, useful for model theoretic applications, and Paulo Agozzini Martin who shows us the beauty of arithmetic and analytic number theory.

\end{document}